\documentclass{amsart}
\usepackage{amsmath,amssymb,amsthm}
\usepackage{ifthen}

\newtheorem{thm}{Theorem}[section]
\newtheorem{lem}{Lemma}[section]
\newtheorem{cor}{Corollary}[section]

\newtheorem{conj}{Conjecture}[section]
\newtheorem{defn}{Definition}[section]

\numberwithin{equation}{section}
\def\Z{\Bbb Z}

\def\R{\Bbb R}

\def\d{\partial}
\def\s{\sigma}
\def\e{\epsilon}

\newcommand{\showcomments}{no}

\newcommand{\commentstyle}{\tiny}

\newcommand{\comment}[1]
{\ifthenelse{\equal{\showcomments}{yes}}
{\footnotemark\marginpar{\sffamily{\commentstyle
\addtocounter{footnote}{-1}\footnotemark#1 }\normalfont}}{}}

\def\flabel#1{\ifmmode #1\else$ #1$\fi}

\let\angle\undefined
\newsymbol\angle 115C
\def\degrees{\ifmmode^\circ\else$^\circ$\fi}
\def\a{\alpha}

\def\g{\gamma}

\def \pict #1 by #2 (#3) {\centerline{
\vbox to #2 {\hrule width #1 height 0pt depth 0pt
\vfill{\special{picture #3 }}}}}

\def \picture #1 by #2 (#3 scaled #4) #5{
\dimen0=#1 \dimen1=#2
\divide\dimen0 by 1000 \multiply\dimen0 by #4
\divide\dimen1 by 1000 \multiply\dimen1 by #4
\vbox{\pict \dimen0 by \dimen1 (#3 scaled #4)
\centerline { #5}}}

\catcode`\@=11

\title[]{The Burnside Ring-Valued Morse Formula for Vector Fields on Manifolds with Boundary }
\begin{document}
\author{Gabriel Katz}
\address{Department of Mathematics, MIT, 77 Massachusetts Ave., Cambridge, MA 02139-4307}
\email{gabkatz@gmail.com}
\date{\today}
\begin{abstract} Let $G$ be a compact Lie group and $A(G)$ its Burnside Ring. For a compact smooth 
$n$-dimensional $G$-manifold $X$  equipped with a generic $G$-invariant vector field $v$, we  prove an equivariant analog of the Morse formula 
\begin{eqnarray}
Ind^G(v) =  \sum_{k = 0}^{n} (-1)^k \chi^G(\d_k^+X)\nonumber
\end{eqnarray}
which takes its values in $A(G)$. Here $Ind^G(v)$ denotes the equivariant index of the field $v$, $\{\d_k^+X\}$ the $v$-induced Morse stratification  (see [M]) of the boundary $\d X$, and $\chi^G(\d_k^+X)$ the class of the $(n - k)$-manifold $\d_k^+X$ \footnote{By definition, $\d_0^+X = X$.} in $A(G)$. \smallskip

We examine some applications of this formula to the equivariant real algebraic fields $v$ in  compact domains $X \subset \R^n$ defined via a generic polynomial inequality.  Next, we link the above formula with the equivariant degrees of certain Gauss maps. This link is an  equivariant generalization of Gottlieb's formulas ([G], [G1]).
\end{abstract}

\maketitle
\section{Introduction}

Let $X$ be a compact smooth $n$-manifold with boundary $\d X$. 
A \emph{generic} (see Definition 1.1) vector field $v$  on $X$ which is nonzero  along $\d X$ gives rise to a  stratification 
\begin{eqnarray}
X := \d_0^+ X \supset \d_1^+X \supset \d_2^+X \supset \dots \supset \d_n^+X 
\end{eqnarray}
by compact submanifolds,  where $dim(\d_j^+X) = n - j$. 
Here  $\d_1^+X$ is the part of the boundary $\d_1X := \d X$ where $v$ points inside $X$. By definition, $\d_2X$ is  the  locus where $v$ is \emph{tangent} to the boundary $\d_1 X$. Its  portion $\d_2^+X \subset \d_2 X$ consists of points  where $v$ points inside $\d_1^+ X$. Similarly, $\d_3X$ is the  locus where $v$ is tangent to $\d_2 X$. In the same spirit, $\d_3^+X \subset \d_3 X$ consists of points where  $v$ points inside  $\d_2^+X$. Continuing this process, we get the Morse stratification (1.1).

In Section 2, Theorem 2.1,  we investigate quite strong restrictions on the nature of the stratification $\{\d_jX\}$ (and thus of $\{\d_j^+X\}$) imposed by a $G$-symmetry of the field $v$ and the manifold $X$. 

\begin{defn}
We say that a field $v$ is \emph{generic}, if for each $k$,  viewed as a section of the bundle $T(\d_kX)|_{\d_{k+1}X}$, $v$ is transversal to the zero section of the tangent bundle $T(\d_{k+1}X)$.
\end{defn}
\smallskip

In his groundbreaking 1929 paper [Mo], Morse discovered a beautiful connection between  stratification  (1.1) and the index $Ind(v)$ of the field $v$. It is expressed in terms of the Euler numbers of the strata from $(1.1)$:
\begin{eqnarray}
Ind(v) = \sum_{k = 0}^n (-1)^k \chi(\d_k^+X). 
\end{eqnarray}
\smallskip

Our main observation is that, for any compact Lie group $G$ and a generic $G$-equivariant  vector field $v$, a similar formula with values in the Burnside ring $A(G)$  holds. Theorem 3.1 is an equivariant version of the Morse formula $(1.2)$ for generic symmetric vector fields on manifolds with boundary. It  is a generalization of Theorem 6.6, [LR], for fields that do not necessarily point outward $X$ along its boundary\footnote{Actually, [LR] deals with the category of cocompact discrete $G$-actions.}. \smallskip

In Section 4, Theorem 4.2, we combine some results of Khovanskii [Kh] about (non-equivariant) indices of real algebraic vector fields in polynomially-defined domains in $\R^n$ with formulas $(3.9)$ and $(3.10)$ from Theorem 3.1 to get a handle on the size of eqivariant indices of such fields.\bigskip

In  Section 5, Theorem 5.1, we obtain an $A(G)$-valued version of  Gottlieb's "Topological Gauss-Bonnet Theorem"  [G], [G1]. These results connect the indices of pullback fields $F^\ast  w$ under  smooth $G$-maps $F: X \to V$ to the equivariant degree $Deg^G(\g)$ of the Gauss map $\g: \d X \to S(V)$. Here $w$ is a  nonsingular $G$-invariant field on a space of a $G$-representation $V$, and  $F|_{\d X}$ is an immersion.
\bigskip

In Section 6, we speculate about  some generalizations of our equivariant index formulas, the generalizations which  reside in refined versions of the ring $A(G)$.
\bigskip

{\it Acknowledgments:\,} I am grateful to Wolfgang L\"uck for a valuable conversation  which launched this investigation.

\section{$G$-invariant Vector Fields and Morse Stratifications}

When $v$ is both generic and $G$-invariant, the Morse stratification $\{\d_k^+X\}$ is invariant as well. In fact, it is quite special. For example, the following proposition is valid:
\begin{thm} Let a compact Lie group $G$ act faithfully on a smooth compact $n$-manifold $X$ with an oriented boundary $\d_1X$. Then, for any $k$ and a generic  $G$-equivariant vector field $v$, 
\begin{itemize}
\item the main orbit-type of $\d_kX$ is $G$. 
\item for all $k > n - dim(G)$, we have $\d_kX = \emptyset$\footnote{That is, $(X, v)$ is $(n + 1 - dim(G))$-convex in the terminology of [K].}. 
\item the sets $\d_{n - dim(G)}^\pm X$ are disjoint unions of \emph{free} $G$-orbits.  
\item if $G$ is connected, the sets $\d_{n - 1 - dim(G)}^\pm X$ are disjoint unions of the following $G$-spaces: 
\begin{enumerate}
\item $G \times [0, 1]$, 
\item $G \times S^1$,
\item the mapping cones $\mathcal C(G \to G/H)$, where $H \approx  SO(2)$ or $SU(2)$,
\item $\mathcal C(G \to G/H) \cup_G\mathcal C(G \to G/K)$, where $H, K \approx  
SO(2)$ or $SU(2)$.
\end{enumerate}
\end{itemize}
\end{thm}

{\it Proof.\;} Pick $x \in \d_kX^\circ := \d_kX \setminus \d_{k+1}X$ and let $H = G_x$. Because $v$ is $H$-invariant, a linearization of the $H$-action at $x$ must leave the flag $F_x$ of  tangent half-spaces $\{T_x^+(\d_l^+ X)\}_{l \leq k}$ invariant as well. Pick a $G$-invariant metric on $X$.  Consider the unique "frame" $\nu_x$ comprised of $k$ mutually orthogonal  rays $r_1, \dots , r_k$  generated by intersecting the flag $F_x$  with the space $N_x^k$, normal in $T_x(X)$ to $T_x(\d_kX)$. Specifically, for all $l < k$,  the positive cone $Span^+\{r_1, \dots ,  r_l  \} = N_x^k \cap T_x^+(\d_{k - l}^+X)$. The "frame" $\nu_x$ must be $H$-invariant. Thus, in a basis consistent with $\nu_x$, the $H$-action is diagonal with the positive eigenvalues.  Since $H$ is compact, its action on the space $N_x^k$ (spanned by $\nu_x$) is trivial. At the same time, the $H$-action on the space tangent to the trajectory $Gx$ at $x$ is trivial as well. \smallskip

Now consider the 0-dimensional $G$-invariant sets $\d_n^\pm X$. By the argument above, for any $x \in \d_n^\pm X$, the group  $H = G_x$ acts trivially on the $n$-flag $F_x$, and thus, on $T_xX$. Since  the $G$-action on $X$ is faithful,  we get $H = 1$. Hence, either  $\d_n^\pm X = \emptyset$  when $dim(G) > 0$, or $dim(G) = 0$ and $\d_n^\pm X$ is a union of free orbits. 

Next, consider the case $dim(G) \geq 1$ and focus on the 1-dimensional $G$-invariant sets $\d_{n - 1}^\pm X$.  By the argument above, for any $x \in \d_{n - 1}^\pm X$, the group  $H = G_x$ acts trivially on the  on the vector space  $N_x^{n -1}$. In addition, $H$ acts trivially   on  the 1-dimensional space tangent to the orbit $G/H$ through $x$. Since $\d_{n -1}X$ is 1-dimensional, the space tangent to the orbit and the space tangent to $\d_{n -1}X$ coincide. As a result, the $H$ action is trivial on $T_xX$, and thus, in the vicinity of $x$. Since  the $G$-action on $X$ is faithful,  we get again $H = 1$. Hence,  either $dim(G) > 1$ and $\d_{n -1}^\pm X = \emptyset$, or $dim(G) =  1$ and $\d_{n -1}^\pm X$ is a union of free $G$-orbits.

Finally, let us treat the general case. For $k > 0$ and a generic $G$-invariant $v$, $\d_kX$ is a closed $G$-manifold. Pick a main orbit $G/H \subset \d_kX$, $H = G_x$ for some $x \in \d_kX$. According to [B], an invariant regular neighborhood of  $G/H \subset \d_kX$ is diffeomorphic to the balanced product $G \times_H V$, where $V$ is a space of an $H$-representation $\psi$.  If this $\psi$ is a non-trivial representation, $G/H$ is not the main orbit:  the stationary groups of some points from $G \times_H (V \setminus \{0\})$ must be smaller than $H$. Thus, we can assume that $H$ acts trivially on $V$. In addition, it acts trivially on the tangent space to the orbit $G/G_x$ at $x$, as well as on the normal space $N_x^k$. As a result, the $G_x$-action on $T_xX$ is trivial. Because $G$ acts faithfully on $X$, $G_x = 1$; that is, the main orbit-type of $\d_kX$ (and hence of $\d_k^\pm X$) is $G$! 
Therefore, when $n - k < dim(G)$, $\d_kX = \emptyset$. 

When $n - k = dim(G)$, each main orbit $G$ is  open and dense in some $G$-invariant union of  connected components of $\d_kX$. Since both $G$ and $\d_kX$ are closed manifolds,  $\d_kX$ is $G$-diffeomorphic to a disjoint union of several $G$'s.

When $n - k = dim(G) + 1$ and $G$ is connected, it is present as a codimension one main orbit-type of each connected component of $\d_kX$. Fortunately, connected closed $G$-manifolds of this kind are very rigid:  they all are either products $G \times S^1$, or unions of two mapping cylinders $\mathcal C(G \to G/H)$, $\mathcal C(G \to G/K)$ with $H$ and $K$ being diffeomorphic to a sphere [B], [K1].  The two mapping cylinders share the same "top" $G$. The only spheres among the Lie groups are $O(1) \approx \Z_2, SO(2), SU(2)$. However, the $\Z_2$-option must be excluded: it leads to a non-orientable manifold, while all $\d_k X$ are oriented for $k > 0$. This leaves us with the models $$\mathcal C(G \to G/H) \cup_G\mathcal C(G \to G/K),$$ where $H, K = SO(2), SU(2)$, for the components of $\d_{n- 1- dim(G)}X$. 

Similarly,  the models for $\d_{n- 1- dim(G)}^\pm X$ are: 
$\mathcal C(G \to G/H) \cup_G\mathcal C(G \to G/K)$, $\mathcal C(G \to G/H)$, $G \times S^1$, and $G \times [0, 1]$, where $H, K = SO(2), SU(2)$. \qed

\begin{cor}Let $X, v$ be as in Theorem 2.1.  Assume that  $G$ is  finite.  Then, for some integer 
$l \geq 0$,\, $|\d_nX| = 2l |G|$, provided  $|G|$ being odd, and $|\d_nX| = l |G|$, provided $|G|$ being even.

If $G = SO(2)$, then $\d_nX = \emptyset$, $\d_{n - 1}^\pm X$ each is a union of circles on which $G$ acts freely, and each $\d_{n - 2}^\pm X$ is a union of tori $SO(2) \times S^1$, cylinders $SO(2) \times [0, 1]$, 2-disks, and 2-spheres. The $SO(2)$-action on each disk or sphere is semi-free. 
\end{cor}

{\it Question.\;} For a given compact $G$-manifold $X$ of dimension $n$ with orientable boundary, what is the minimal number of free trajectories that form the sets $\d_{n - dim(G)}(X)$ and $\d_{n - dim(G)}^\pm(X)$? The minimum is taken over the space of all generic $G$-invariant vector fields $v$\footnote{In [K] we proved that, for $n = 3$ and $G =1$, this minimum is zero.}.

\begin{lem} For a $G$-invariant field $v$ and each compact subgroup $H \subset G$, we have
$\d_k^+(X^H, v|_{X^H}) = X^H \cap \d_k^+(X, v)$.
\end{lem}

{\it Proof.\;} 
A key observation here is that, for each compact subgroup $H \subset G$, any $G$-invariant field $v$ on $X$ is tangent to the submanifold $X^H \subset X$. Indeed, in a $G$-invariant metric $g$, the normal to $X^H$ component $\nu_H$ of $v$ must be $H$-invariant, and thus vanishes for $H \neq 1$. 
Each submanifold $X^H$ meets $\d_1X$ transversally: just consider a linearization of the $G_x$-action, $G_x \supset H$,  at a typical point $x \in X^H \cap \d_1X$. Since $v$ is tangent to $X^H$,  evidently, $\d_1^+(X^H, v|_{X^H}) = X^H \cap \d_1^+(X, v)$. 

We proceed by induction on $k$. Assume that $X^H$ and $\d_sX$ are transversal,  $\d_s(X^H, v|_{X^H}) = X^H \cap \d_s(X, v)$, and $\d_s^+(X^H, v|_{X^H}) = X^H \cap \d_s^+(X, v)$ for all $s < k$. If $x \in X^H \cap \d_k^+(X, v)$ and  $v(x)$  is tangent to $\d_{k -1}X \cap X^H = \d_{k-1}X^H$, then $v(x)$ points inside of $\d_{k -1}^+X$ if and only if it points inside $\d_{k -1}^+X^H$. Thus, $\d_k(X^H, v|_{X^H}) = X^H \cap \d_k(X, v)$, which completes the induction step. \qed
\bigskip

\section{Equivariant Morse Formula for Vector Fields}

Following tom Dieck, consider  a  ring $A(G)$ generated over $\Z$ by  equivalence classes of  compact
differentiable $G$-manifolds. Two $G$-manifolds $X$ and $Y$ are said to be equivalent, if for any compact subgroup $H \subset G$, the Euler numbers $\chi(X^H)$ and $\chi(Y^H)$ are equal. This equivalence relation respects disjoint unions and cartesian products of $G$-spaces. Therefore, $A(G)$ is a ring with the sum and product operations induced by disjoint unions and cartesian product of equivalence classes of $G$-manifolds. 

In particular, any compact $G$-manifold $X$ gives rise to an element $\chi^G(X) \in A(G)$.

Let  $ch: A(G)\to \prod_{(H) \in conj(G)} \; \Z$ be a ring homomorphism induced by the correspondence $\chi^G(X) \to   \{\chi(X^H)\}_{(H) \in conj(G)}$, where  $conj(G)$ denotes the conjugacy classes of compact subgroups in $G$. It turns out that $ch$ is a \emph{monomorphism} [D]. We denote by $ch_H$ the $(H)$-indexed component of $ch$.
\smallskip

Consider only the conjugacy classes $(H)$ with the finite quotients $WH :=  
NH/H$ ($NH$ standing for the normalizer of $H$ in $G$) and denote by $\Phi(G)$ the set of such $(H)$. 
Next, form a free abelian group $A'(G)$ generated by elements of $\Phi(G)$ (equivalently, by the orbit-types $\{G/H\}_{(H) \in \Phi(G)}$ regarded as $G$-spaces). By [D], Theorem 1, the natural  homomorphism $A'(G) \to A(G)$ is an isomorphism. As a result, any element $[X] \in A(G)$ is detected by 
$\{ch_H([X])\}_{(H) \in \Phi(G)}$ alone.
\bigskip

Let $U_Z$ be a $G$-invariant regular neighborhood\footnote{In particular, $Z$ is an equivariant deformation retract of $U_Z$} of the zero set $$Z := Z(v) = \{x \in X|\; v(x)=0\}$$ and a smooth manifold. Assume that $v$ is in general position with respect to $\d_1U_Z := \d U_Z$. Note that, for each $H \subset G$, $U_Z^H$ and $Z^H$ are homotopy equivalent. \smallskip

Inspired by [M], we propose the following
\begin{defn}
\begin{eqnarray}
ind^G(v) := \sum_{k = 0}^{n} (-1)^k \chi^G(\d_k^+U_Z) = \chi^G(Z) + \sum_{k = 1}^{n} (-1)^k \chi^G(\d_k^+U_Z)
\end{eqnarray}
\end{defn}

\begin{lem} For each $(H) \in conj(G)$, $ch_H(ind^G(v)) =  ind(v|_{X^H})$, the non-eqivariant index of the field $v$ on $X^H$ (equivalently, on $U_Z^H$).
\end{lem}
{\it Proof.\;} By Lemma 2.1, $v$ is tangent to $X^H \supset U_Z^H$ and $(\d_k^+U_Z)^H = \d_k^+(U_Z^H)$. Since $v$ is tangent to  $X^H$,  all the zeros of $v|_{X^H}$ are among the zeros of $v$ on $X$, and no new zeros appear. By the Morse formula $(1.2)$ ([M]), $$ind(v|_{U_Z^H}) = \sum_{k = 0}^{n} (-1)^k \chi(\d_k^+(U_Z^H)) = 
\sum_{k = 0}^{n} (-1)^k \chi(\d_k^+(U_Z)^H).$$
Thus, 
\begin{eqnarray}
ch_H(ind^G(v)) =  \sum_{k = 0}^{n} (-1)^k ch_H((\d_k^+U_Z)) = ind(v|_{U_Z^H}) = ind(v|_{X^H}). 
\end{eqnarray}
\qed
\bigskip

For $G$-invariant vector fields, there is an alternative definition of the equivariant index.  It is more involved and mimics the classical definition of index as the sum of local degrees produced by the field $v$ in the vicinity of its zero set.  Let us describe this definition along the lines of [L], [LR].\smallskip

The tangent space $T_{(x, 0)}(TX)$ decomposes as a direct sum of the subspaces $T_xX \oplus T_xX$, where the first factor is thought as being tangent to the zero section $\e: X \to TX$ and the second one  to the fiber of the bundle $TX \to X$. Let $ \a: T_xX \oplus T_xX \to T_{(x, 0)}(TX)$ denote this $G_x$-equivariant isomorphism. Consider a $G$-equivariant  vector field $v$ on $X$, viewed as a section $v: X \to TX$ of the tangent $G$-bundle $TX$.  Assume that $v$ is transversal to the zero section $\e$ (due to the "doubling" nature of  the $G_x$-space $T_{(x, 0)}(TX)$, this equivariant transversality of $v$ and $\e$ is available by a general position argument). If $v(x) = 0$ for some $x \in X$, then  $T_{(x, 0)}(TX)$ decomposes a direct sum of $Dv(T_xX)$ and $D\e(T_xX)$ as well.  Here $Dv$ and $D\e$ denote the differentials of the respective maps. Consider the $G_x$-equivariant isomorphism 
\begin{eqnarray}
\mathcal D_{x, v}: T_xX \stackrel{Dv}{\longrightarrow} T_{(x, 0)}(TX)  \stackrel{\a^{-1}}{\longrightarrow}  T_xX \oplus T_xX  \stackrel{pr_2}{\longrightarrow}  T_xX
\end{eqnarray}
, where $pr_2$ stands for the projection on the second factor. $\mathcal D_{x, v}$ induces a $G_x$-map $\mathcal D_{x, v}^c: T_xX^c \to T_xX^c$ on the one point compactification $T_xX^c$ of $T_xX$. For such a map, its degree,  $Deg^{G_x}(\mathcal D_{x, v}^c) \in A'(G_x) \approx A(G_x)$, is defined in terms of the equivariant Lefschetz classes $\Lambda^{G_x}(\mathcal D_{x, v}^c)$,  
$\Lambda^{G_x}(Id_{T_xX^c})$ by the formula 
\begin{eqnarray}
Deg^{G_x}( \mathcal D_{x, v}^c) =  \big [\Lambda^{G_x}(\mathcal D_{x, v}^c) - 1\big ] \big [\Lambda^{G_x}(Id_{T_xX^c}) - 1\big ]
\end{eqnarray}
(see [L],  [LR]), where $1 \in A(G_x)$ stands for the class of a point.\smallskip 

In Section 4, we will return briefly to the definition of an equivariant degree $Deg^G(f)$ for a general $G$-map $f: X \to Y$ of two $G$-manifolds.\smallskip

Meanwhile, let us recall the definition of  the equivariant Lefschetz class $\Lambda^{G}(f) \in A'(G)$ of a $G$-map $f: X \to X$, $X$ being a finite $G-CW$-complex\footnote{for example, see [O] for the definition of a $G$-$CW$-complex}. In particular, for each $(H) \in \Phi(G)$, the orbit-space $WH\setminus X^H$ has the structure of a finite $CW$-complex. 
Let
\begin{eqnarray}
\Lambda^G(f)  := \sum_{(H) \in \Phi(G)} \; L^{\Z [WH]}(f^H, f^{> H}) \cdot \chi^G(G/H)
\end{eqnarray}
, where $\chi^G(G/H)$ is the class of the homogeneous space $G/H$ in $A'(G)$,  and the integer 
\begin{eqnarray}
L^{\Z [WH]}(f^H, f^{> H}) = \sum_{p \geq 0} \; (-1)^p \cdot tr_{\Z [WH]}\big(C_p\big(f^H, f^{> H})\big).
\end{eqnarray}
In $(3.6)$, we employ the $\Z$-valued trace $tr_{\Z [WH]}$ of the $f$-induced automorphism \hfill\break
$C_p\big(f^H, f^{> H})$ of the finitely generated free 
$\Z [WH]$-module $C_p\big(X^H, X^{> H}; \Z)$, the module of the relative celluar $p$-chains of the pair $(X^H, X^{> H})$ (see [LR] for details).
\smallskip

\begin{lem}
The equivariant Lefschetz class $\Lambda^{G}(f)$ from $(3.5)$ is detected by the non-equivariant Lefschetz numbers $\{\Lambda(f^H) \in \Z\}_{(H) \in \Phi(G)}$. 
\end{lem}

{\it Proof. \;} By [D], Theorem 2,  only the orbit-types $G/H$ with $WH$ being finite contribute to the homomorphism $ch:  A(G) \to \prod_{(H)} \; \Z$.  Therefore, we will consider only the orbit-types from $\Phi(G)$. \smallskip

In (3.5), we took into account only the orbit-types $(H) \in \Phi(G)$ with non-vanishing coefficients $L^{\Z [WH]}(f^H, f^{> H})$ (defined by $(3.6)$). Among them, pick a minimal\footnote{i.e. a \emph{maximal} subgroup $K \subset G$}  orbit-type $(K)$ and apply $ch_K$ to $\Lambda^{G}(f)$. Then,

\begin{eqnarray}
ch_K(\Lambda^G(f)) = \sum_{(H) \in \Phi(G)}\Big\{\sum_{p \geq 0} \; (-1)^p\; tr_{\Z[WH]}\big(C_p(f^H, f^{> H})\big) \Big\} \cdot \chi(G/H^K) \nonumber \\
= \Big\{\sum_{p \geq 0} \; (-1)^p\; tr_{\Z[WK]}\big(C_p(f^K, f^{> K})\big)   \Big\} \cdot \chi(G/K^K) =  \nonumber \\
= \Big\{\sum_{p \geq 0} \; (-1)^p\; tr_{\Z[WK]}\big(C_p(f^K, f^{> K})\big)   \Big\} \cdot |WK| =  \nonumber \\ 
= \sum_{p \geq 0} \; 
(-1)^p\; tr_{\Z}\big(C_p(f^K, f^{> K})\big) := \Lambda(f^K).  
\end{eqnarray}
In other words, for such a minimal $(K) \in \Phi(G)$, $$L^{\Z[WK]}(f^K, f^{> K}) = |WK|^{-1} \Lambda(f^K) = |WK|^{-1} ch_K(\Lambda^G(f)).$$

The rest of the argument is performed inductively, the induction step being applied to a minimal orbit-type $(K')$  in $\Phi(G) \setminus (K)$ with a non-zero coefficient $L^{\Z [WK']}(f^{K'}, f^{> K'})$. Indeed,  compute $ch_{K'}$ of $$\beta := \Lambda^G(f) - L^{\Z[WK]}(f^K, f^{> K})\cdot \chi^G(G/K)$$ to conclude, in a similar way,   that $ch_{K'}(\beta)$ equals $\Lambda(f^{K'})$.
\qed
\smallskip

\begin{cor}  $Deg^{G_x}(\mathcal D_{x, v}^c) \in A(G_x)$ is determined by the $\Z$-valued degrees \hfill\break $\big\{\,Deg(\mathcal D_{x, v}^{H, c})\,\big\}_{(H) \in \Phi(G_x)}$  of the maps $\mathcal D_{x, v}^{H, c} := \mathcal D_{x, v}^c| : (T_xX^c)^H \to (T_xX^c)^H$. 
For all $H \notin \Phi(G_x), \; H \subset G_x$,  $ch_H [Deg^{G_x}(\mathcal D_{x, v}^c)] = 0$.
\end{cor}

{\it Proof.\;} Apply the  arguments (3.7) of Lemma 3.2 (with $G$ being replaced by $G_x$) to formula $(3.4)$.
\qed
\bigskip

Now put
\begin{eqnarray}
Ind^G(v) := \sum_{Gx \in \,G\setminus Z(v)} \;  G \times_{G_x} [Deg^{G_x}(\mathcal D_{x, v}^c)]
\end{eqnarray}

\begin{lem} $Ind^G(v)$, defined by $(3.8)$, and $ind^G(v)$, defined by $(3.1)$, are equal.
\end{lem}
{\it Proof.\;} 
For each $x$ from the orbit-space  $G \setminus  Z(v)$, consider a small $G_x$-equivariant  $n$-disk  $D_x \subset X$ centered on $x$. Note that, when $Z(v)$ is discrete, the quotient  $G/G_x$ must be finite. Evidently, $\d_k^+U_Z \approx \coprod_x G \times_{G_x} (\d_k^+D_x)$. Therefore, in order to prove $Ind^G(v) = ind^G(v)$, it will suffice to verify that, for each $x \in G \setminus Z(v)$, $$Deg^{G_x}(\mathcal D_{x, v}^c) = \sum_k (-1)^k \chi^{G_x}(\d_k^+D_x)$$ in $A'(G_x)$. 
The last equality follows from Corollary  3.1 together with the validity of non-equvariant Morse formulas of the type $(1.2)$ (with $G$ being replaced by $H \subset G_x$ and $X$ by $D_x$).
\qed

\begin{thm} Let $G$ be a compact Lie group, and let  $X$ \footnote{also denoted by $\d_0^+X$} be a compact smooth $n$-dimensional $G$-manifold with orientable boundary  $\d_1X$. Let $v$ be a generic $G$-invariant vector field on $X$ (with isolated singularities). Then the equivariant Morse formula 
\begin{eqnarray}
Ind^G(v) = ind^G(v) = \sum_{k = 0}^{n} (-1)^k \chi^G(\d_k^+X)
\end{eqnarray}
is valid in the ring $A(G)$.
If the $G$-action on $X$ is faithful and $dim(G) > 0$, 
\begin{eqnarray}
Ind^G(v) = ind^G(v) = \sum_{k = 0}^{n - 1- dim(G)} (-1)^k \chi^G(\d_k^+X)
\end{eqnarray}
When $G$ is connected, the contribution of the stratum $\d_{n - 1- dim(G)}^+X$ to $Ind^G(v)$ is quite special:
$$\chi^G(\d_{n - 1- dim(G)}^+X) =  \sum_{\{(H) \in \Phi(G) |\, H  \approx SO(2), SU(2)\}} \; 
n_H\cdot \chi^G(G/H),$$
where $\{n_H\}$ are some nonnegative integers. In particular, if a connected $G$ is such that,  for each  $H \subset G$ isomorphic to $SO(2)$ or $SU(2)$, the group $WH$ is infinite, 
then the contribution of the stratum $\d_{n - 1- dim(G)}^+X$ vanishes. 
\end{thm} 

{\it Proof.\;} By Lemma 3.3, $Ind^G(v) = ind^G(v)$. Since the elements of $A(G)$ are detected by the character map $ch: A(G) \to \prod_{(H) \in \Phi(G)} \; \Z$,  it will suffice to check the validity of (3.9) by applying $ch$ to both sides of the conjectured equation $ind^G(v) = \sum_{k = 0}^{n} (-1)^k \chi^G(\d_k^+X)$ . 

By Lemma 2.1, for any $H \subset G$, $ch_H[\chi^G(\d_k^+X)] := \chi((\d_k^+X)^H) = \chi(\d_k^+(X^H))$. Thus, 
\begin{eqnarray}
ch_H\big(\sum_{k = 0}^{n} (-1)^k \chi^G(\d_k^+X)\big) = \sum_{k = 0}^{n} (-1)^k \chi(\d_k^+(X^H)) 
\end{eqnarray}

On the other hand, by Lemma 3.1 and formula $(3.2)$, $ch_H(ind^G(v)) =  ind(v|_{X^H})$. Since the Morse formula claims that $ind(v|_{X^H}) = \sum_{k = 0}^{n} (-1)^k \chi(\d_k^+(X^H))$, we get 
the desired equality $ch_H(ind^G(v)) = ch_H(\sum_{k = 0}^{n} (-1)^k \chi^G(\d_k^+X))$ for all  $(H) \in conj(G)$. \smallskip

In order to derive formula $(3.10)$ from formula $(3.9)$, we employ Theorem 2.1 to conclude that $\d_k^+X = \emptyset$ for all $k > n  - dim(G)$. Moreover, $\d_{n - dim(G)}^+X$ is the union of free $G$-orbits, and thus $\chi^G(\d_{n - dim(G)}^+X) = 0$, provided $dim(G) > 0$. 

The last claim of the theorem can be validated by the following observation. Models $(1)$ and $(2)$ from Theorem 2.1 have zero  classes in $A(G)$, model (3) is $G$-homotopy equivalent to $G/H$, and model (4), by the additivity of Euler characteristics,  produces the element $\chi^G(G/H) + \chi^G(G/K)$. 

If a connected $G$ is such that,  for each  $H \subset G$ isomorphic 
$SO(2)$ or $SU(2)$, $WH$ is infinite, the space $G/H$ admits a free $S^1$ action, $S^1 \subset WH$. Thus, $\chi^G(G/H) =0$. As a result, for  such a $G$, we get a simplification of $(3.10)$:
\begin{eqnarray}
Ind^G(v)  = \sum_{k = 0}^{n - 2- dim(G)} (-1)^k \chi^G(\d_k^+X).
\end{eqnarray}
\qed
\bigskip

\begin{cor}Let $X$ and $v$ be as in Theorem 3.1. Denote by  $v_1^+$  an orthogonal projection (with respect to a $G$-invariant metric on $X$) of the field $v|_{\d_1^+X}$ on the tangent space $T(\d_1^+X)$. Assume that $v_1^+$ has only isolated singularities. Then the following formula holds in $A(G)$:
\begin{eqnarray}
\chi^G(X) = ind^G(v) + ind^G(v_1^+) 
\end{eqnarray}
\end{cor}

{\it Proof.\;} In view of the recursive nature of the Morse stratification $\{\d_k^+X \}$, $(3.9)$ and $(3.13)$, applied to all consecutive terms in the stratification, are equivalent formulas.
\qed

\section{Polynomial Vector Fields in Polynomial Domains}

Now, we turn our attention to polynomial vector fields $v$ in domains  in $\R^n$ that are defined by polynomial inequalities.\smallskip

Consider a real $n$-dimensional vector space $W$ and a $G$-representation $\Psi: G \to GL_{\R}(W)$. Denote by $\mathcal P_\Psi$ the algebra of invariant polynomials on $W$.  Let $v$ be a $G$-invariant vector field in $W \approx \R^n$  which has polynomial components $\{P_i \in \mathcal P_\Psi\}_{1 \leq i \leq n}$.  Assume that  $deg(P_i) \leq m_i$. We denote by $\tilde P_i$ the homogenized versions of $P_i$, that is, $\tilde P_i(1, x_1, \dots , x_n) = P_i(x_1, \dots , x_n)$.

Also consider an invariant polynomial $Q \in \mathcal P_\Psi$ of degree $d$. We say that the pair $(v, Q)$ is \emph{non-degenerated}, if 1) all zeros of $v$ are simple,  2) the system $\{x_0 = 0, \; \tilde P_i(x_0, x_1, \dots, x_n) = 0\}_{1\leq i \leq n}$ has only a trivial solution\footnote{i.e. no zero of $v$ escapes to infinity.}, and 3) the hypersurface $Q = 0$ does not contain the zeros of $v$. 

For any non-increasing sequence $m_1, m_2, \dots , m_n$ of natural numbers, form the parallelepiped $\Pi(m_1, \dots, m_n)$ in $\R^n$ defined by
\begin{eqnarray}
\{0 \leq x_i \leq m_i - 1\}_{1\leq i \leq n}
\end{eqnarray}

Let $O(d, m_1, \dots,  m_n)$ be the number of integral lattice points $(x_1, \dots , x_n)$  in  $\Pi(m_1, \dots m_n)$,  subject to the inequalities 
\begin{eqnarray}
\qquad  \frac{1}{2}(m_1 + \dots + m_n - d - n) \leq x_1 + \dots + x_n \leq \frac{1}{2}(m_1 + \dots + m_n  - n)
\end{eqnarray}

Consider the domain $X_Q := \{w \in W| \; Q(w) \geq 0\}$. Then, according to a theorem of Khovanskii  (Theorem 1, [Kh]), for a non-degenerated pair $(v, Q)$ as above, the absolute value of index, $| ind(v) |$, in $X_Q$ is bounded from above by  $O(d, m_1, \dots,  m_n)$; moreover,  this estimate is sharp. 

Assume that $X_Q$ is compact with a smooth boundary $\d_1X_Q$ and a polynomial field $v$ is generic (see Definition 1.1) in relation to the boundary. Combining the Khovanskii Theorem and Morse Formula $(1.2)$, we get

\begin{thm} For a polynomial field $v$ in $X_Q$ as above,
\begin{eqnarray}
\Big|\sum_{k = 0}^n \; (-1)^k \chi(\d_k^+X_Q) \Big| \leq O(d, m_1, \dots,  m_n)
\end{eqnarray}
Moreover, there exist $(Q, v)$ for which the inequality $(4.3)$ can be replaced by the equality.
\end{thm}

Let $V \subset \R^n$ be a vector subspace of dimension $l$. Let $d_V \leq d$ denote the degree of the polynomial $Q$ being restricted to the subspace $V$. Consider a subspace $U \subset \R^n$ which is spanned by some set of $l$ basic vectors $\{e_{j_1}, \dots e_{j_l}\}$ in $\R^n$ and such that the obvious orthogonal projection $p_U: \R^n \to U$, being restricted to $V$, is onto. Denote by $\mathcal U(V)$ the finite set of such $U$'s.

Put 
\begin{eqnarray}
O(V; d_V, m_1, \dots , m_n) = min_ {U \in \mathcal U(V)} \Big\{O(d_V, m_{j_1}, \dots , m_{j_l})\Big\},
\end{eqnarray}
 where $U = span\{e_{j_1}, \dots, e_{j_l} \}$ and $O(d_V, m_{j_1}, \dots , m_{j_l})$ is defined as in $(4.2)$.
\bigskip

The theorem below  generalizes the estimates $(4.3)$ for an equivariant setting.
\begin{thm}
Let $G$ be a compact Lie group. Pick a sequence $m_1 \geq m_2 \geq \dots \geq m_n > 0$ of integers. Consider a representation $\Psi: G \to GL(n, \R)$,  an invariant polynomial $Q(x_1, \dots , x_n)$ of degree $d$, and  a non-degenerate 
$\Psi(G)$-invariant polynomial field $$v = (P_1(x_1, \dots , x_n),\, \dots ,\, P_n(x_1, \dots , x_n))$$ such that $deg(P_i) \leq m_i$.  Assume that $X_Q := \{\vec x \in \R^n | \; Q(\vec x) \geq 0\}$ is a compact\footnote{For example, take $Q = a_1 x_1^d + \dots + a_n x_n^d \in \mathcal P_\Psi$ (all $a_i >  0$)  plus a lower degree polynomial.} domain with a smooth boundary $\d_1X_Q$ and that $v$ is in general position with respect to $\d_1X_Q$. 

Then the image $\prod_{(H) \in \Phi(G)} z_{(H)}$ of the element $\sum_{k = 0}^n \; (-1)^k \chi(\d_k^+X_Q) \in A(G)$ under the character monomorphism $ch: A(G) \to \prod_{(H) \in \Phi(G)} \Z$ belongs to the parallelepiped $\mathcal P_{ \Phi(G)}$ defined by the inequalities:
\begin{eqnarray}
\big | z_{(H)} \big | \leq O((\R^n)^H; d_{(\R^n)^H}, m_1, \dots,  m_n) 
\end{eqnarray}
\end{thm}

{\it Proof.\;} Let $l = dim((\R^n)^H)$. For each $U = span\{e_{j_1}, \dots, e_{j_l} \}$  as above with the property $p_U: (\R^n)^H \to U$ being onto, the projection $p_U$ induces an invertible linear transformation of the pairs 
$(v,\; X_Q ^H)$ and $(p_U(v),\; p_U(X_Q^H))$, where $X_Q^H := X_Q \cap (\R^n)^H$. Let $I$ be the set of indices complementary to the set $J := \{j_1, \dots , j_l\}$. Then the components $\tilde P_j$ of $p_U(v)$ are obtained  by substituting $\{x_i = 0\}_{i \in I}$ into $\{P_j(x_1, \dots , x_n)\}_{j \in J}$. Note that $deg(\tilde P_j) \leq deg(P_j) \leq m_j$. Also, $p_U(X_Q^H) \subset Span\{e_j\}_{j \in J}$ is defined by a polynomial inequality $\tilde Q(x_{j_1}, \dots , x_{j_l}) \geq 0$ which is obtained from $Q(x_1, \dots , x_n) \geq 0$ by a substitution that expresses each $x_i,\, i \in I,$ as linear combination of the $\{x_j\}_{j \in J}$. Evidently, $deg(\tilde Q) = d_{(\R^n)^H}  \leq  d$. Since $v|_{(\R^n)^H}$ is parallel to $(\R^n)^H$,  $p_U$ also maps $\d_k^+(X_Q^H)$ onto $\d_k^+[p_U(X_Q^H)]$. By [Kh] and Theorem 3.1, $|Ind(p_U(v))|$ in $p_U(X_Q^H)$ (equivalently, $|\sum_{k = 0}^n \; (-1)^k \chi(\d_k^+p_U(X_Q))|$) has an upper boundary $O(d_{(\R^n)^H}, m_{j_1}, \dots , m_{j_l})$. Via the linear diffeomorphism $p_U$, $Ind(p_U(v))$ in $\{\tilde Q \geq 0\}$  and $Ind(v|_{X_Q^H})$ in $X_Q^H$ are equal. Thus, $|Ind(p_U(v))| \leq O(d_{(\R^n)^H}, m_{j_1}, \dots , m_{j_l})$. In view of formula $(4.3)$ and Lemma 3.1, $$|ch_H(Ind^G(v))| =  |Ind(v|_{X_Q^H})| \leq O(d_{(\R^n)^H}, m_{j_1}, \dots , m_{j_l}),$$ and thus using formula-definition $(4.4)$, $$|ch_H(Ind^G(v))| \leq O((\R^n)^H; d_{(\R^n)^H}, m_1, \dots , m_n).\qed$$

\begin{cor} Let $v$ be as in Theorem 4.2 and  $\Psi$ be an orthogonal representation. Denote by $B_r \subset \R^n$ the ball of radius $r$ centered at the origin. Put $d_H= dim((\R^n)^H)$. Then, for each $(H) \in \Phi(G)$ and generic $r$, $$\big|ch_H(Ind^G(v|_{B_r}))\big|  = \big|\sum_{k = 0}^{d_H} \; (-1)^k \chi(\d_k^+B_r^H)\big| \leq O((\R^n)^H; 2 , m_1, \dots , m_n).$$
\end{cor}
{\it Proof.\;} Take $x_1^2 + \dots + x_n^2 - r^2$ for the role of $Q$ from Theorem 4.2 and pick $r$ so that $v$ is generic with respect to $\d_1B_r$. 
\qed

\section{The Gottlieb  and  Gauss-Bonnet Equivariant Formulas}

Our next goal is to reinterpret Gottlieb's formulas ([G1]) for indices of pullback vector fields within an equivariant setting.  \smallskip

Let us recall the notion of a \emph{pullback} field (see [G1]). Let $F: X \to Y$ be a differentiable map of two $n$-dimensional Riemannian manifolds, and $w$ a vector field on $Y$. Let $F^\ast w$ be a field on $X$ defined by the formula
$$\big< (F^\ast w)(x), \; u(x) \big>_X \; = \; \big<w, \; DF(u(x)) \big>_Y$$
, where $x \in X$,  $u(x) \in T_xX$ is a generic vector, and $\big < \sim, \sim \big >$ denotes the scalar product in the appropriate tangent space. In other words, if $\omega$ is a 1-form dual  to $w$ in $Y$, then $F^\ast w$ is dual to $F^\ast\omega$ in $X$. 

In the case of a gradient field $w = \nabla f$ on $Y$ ($f: Y \to \R$ being a smooth function), 
$F^\ast w = \nabla(f \circ F)$.\smallskip

Note that when $F$ is an equivariant map,  both metrics on $X$ and $Y$ are $G$-invariant, and $w$ is an invariant field, then $F^\ast w$ is invariant as well.\smallskip

We recall the notion of an equivariant degree $Deg^G(F) \in A(G)$ of a $G$-map $F: X \to Y$ between two compact $G$-manifolds of the same dimension (cf. [L]). Crudely, it is an element of $\prod_{(H) \in conj(G)} \; \Z$ whose $(H)$-component is the usual $deg(F^H)$, where $F^H: X^H \to Y^H$. In fact, such an element $Deg^G(F) \in A(G)$. This naive construction runs into some complications because of the ambiguities in choosing orientations of fixed point components, both in the source and the target. In a sense, one wants to coordinate the orientations of $X^H$ and $Y^H$ (when they are of the same dimension). The issues with the coherent orientations can be resolved by introducing some additional synchronizing structure called in [L] "the $O(G)$-transformation of the fiber transports". Roughly speaking, it assigns a transfer $G$-map $(T_{F(x)}Y)^c \to (T_xX)^c$ to each $x \in X$. \smallskip

Fortunately, we need to  employ $Deg^G(F)$ in a particular situation, where the synchronization of the orientations can be achieved by pedestrian means.  Consider an equivariant immersion $f$ of a closed oriented $(n - 1)$-dimensional $G$-manifold $Z$ into an real $n$-dimensional  space $V$ of an orthogonal $G$-representation. Denote by $S(V)$ the unit sphere centered at the origin. To  each 
$x \in Z$ we assign the unit vector $n(x)$ tangent to $V$ at $f(x)$ and normal to the $f$-image of a  small neighborhood $U_x \subset Z$ of $x$. 
The orientations of $V$ and $Z$ help to resolve the two-fold ambiguity in picking $n(x)$.

Since $f$ is an immersion, for any $H \subset G$, we get $Z^H = f^{-1}(f(Z) \cap V^H)$. Moreover, because $n(x)$ is orthogonal to $f(U_x)$ at $f(x) \in V^H$, $n(x)$ must be parallel to $V^H$. Indeed, if $n(x)$ would have a nontrivial component $\nu_x$ normal to $V^H$, $\nu_x$ must be moved by elements of $H \setminus \{1\}$; on the other hand, $n(x)$ is $H$-invariant since $T_{f(x)}(U_x)$ is. 

In fact, $f^H: Z^H \to V^H$ is an immersion as well. Therefore an orientation of $V^H$, with the help of the "field" $n(x)$ (we assume that $x \in Z^H \setminus Sing(f^H)$), picks a particular orientation  of $Z^H$. In the following, we assume that the orientations of $V^H$ and $Z^H$ are always synchronized in this way. 
\smallskip

Thus  the Gaussian map $\g^H: Z^H \to S(V^H)$ is well-defined for any $(H)$ that occurs as an orbit-type of $Z$. Unless $dim(V^H) \leq 1$, $S(V^H)$ is connected.  In such case, the degree $deg(\g^H)$ is defined as a sum of degrees of maps $\{\g^H_\a: Z^H_\a \to S(V^H)\}_\a$, where $Z^H_\a$ denotes a typical connected component of $Z^H$. When $dim(V^H) = 1$,  $S(V^H) = a \coprod b$, and  $deg(\g^H)$ is defined as the sum of  degrees of the two obvious maps with the singleton targets $a$ and $b$.  
\smallskip

The theorem below is an equivariant version of the "Topological Gauss-Bonnet Theorem" from [G, page 466].
\begin{thm} 
Let $V$ be a real vector space of dimension $n$ on which a compact Lie group $G$ acts orthogonally. We assume that $V$ admits an invariant non-vanishing vector field $w$\footnote{For example, such $w \neq 0$ exists when $V = U \oplus \R^1$, where $U$ is the space of a $G$-representation: just put $w = \nabla(f)$, where $f: U \oplus \R^1 \to \R^1$ is the obvious projection.}.  Let $X$ a compact smooth $n$-dimensional $G$-manifold with an oriented boundary $\d_1X$.  Consider a $G$-map  $F: X \to V$ whose Jacobian is non-zero on $\d_1X$. 
Let $v = F^\ast w$ be the the pullback of $w$ under $F$. Denote by  $\{\d_k^+X\}_{0 \leq k \leq n}$ the Morse stratification of $X$ induced by $v$.

Then the degree of the Gauss map $\g: \d_1X \to S(V)$ with values in  $A(G)$ can be computed by
\begin{eqnarray}
Deg^G(\g) = \chi^G(X) - Ind^G(v) = - \sum_{k = 1}^{n} (-1)^{k} \chi^G(\d_k^+X)
\end{eqnarray}
Hence, $Ind^G(v)$ and the RHS of $(5.1)$ are $w$-independent\footnote{Note that $(5.1)$ proves that $Deg^G(\g)$, a priori an element of $\prod_{(H) \in conj(G)} \; \Z$, actually  belongs to $A(G)$.}.
\end{thm} 
{\it Proof.\;} Since $DF|_{\d_1X}$ is of the maximal rank, we can pullback the $G$-invariant riemannian metric $g$ in $V$ to an equivariant collar of $\d_1X \subset X$ and then extend the pullback $F^\ast(g)$ to an invariant metric on $X$. Let $n$ be the unitary field outward normal to $\d_1X$. Then, as we described prior to the statements of Theorem 5.1, the Gauss map $\g: x \to DF(n(x))$,  $x \in \d_1X$,  is well-defined, equivariant, and helps to pick coherent orientations of components in $\d_1(X^H)$. We have noticed already that $n(x)$ is contained in $T_x(X^H)$ and is normal to $\d_1(X^H)$, provided $x \in X^H$. Also, for $x \in X^H$,  $v(x) \in T_x(X^H)$. We can apply the non-equivariant Gottlieb's  formula to each $\g^H: \d_1(X^H) \to S(V^H)$ to conclude that $deg(\g^H) = \chi(X^H) - Ind(v|_{X^H}) = - \sum_{k = 1}^{dim(V^H)} (-1)^{k} \chi(\d_k^+(X^H))$. By Lemma 2.1, the latter sum is the $(H)$-component $ch_H$ of  $- \sum_{k = 1}^{n} (-1)^{k} \chi^G(\d_k^+(X)) \in A(G)$. Finally, $\{deg(\g^H)\}_{(H) \in \Phi(G)}$ detect $Deg^G(\g)$. Along the way, we have shown that $Deg^G(\g) \in A(G)$. \qed
\bigskip

{\it Remark\;} Formula $(5.1)$ tells us that we can \emph{define} $Deg^G(\g) \in A(G)$ as 
\begin{eqnarray}
- \sum_{k = 1}^{n} (-1)^{k} \chi^G(\d_k^+X)
\end{eqnarray} 
for any choice of an invariant field $w \neq 0$ in $V$, thus a priori avoiding all the troubles with the orientations. On the other hand, this definition of $Deg^G(\g)$ makes perfect sense for any equivariant immersion $f: Z \to V$ of a closed oriented $G$-manifold $Z$ of codimension one: just define $\d_0^+Z$ to be the locus  $\{x \in Z|\; \big <n(x), w(x)\big > \leq 0\}$, and then proceed as in $(1.1)$.  We conjecture that, in general,  
$\sum_{k = 0}^{n - 1} (-1)^{k} \chi^G(\d_k^+Z)$ is $w$-independent. To prove this conjecture will suffice to construct an equivariant coboundary $X$  for $Z$ and to extend $f: Z \to V$ to a $G$-map $F: X \to V$.

\begin{cor} Let $g$ denotes the $G$-invariant Riemannian metric on $V$. Under notations and hypotheses of Theorem 5.1 and  we get
\begin{eqnarray}
\int_{\d_1X^H} K_H\, d\mu_H = - \sum_{k = 1}^{dim(V^H)} (-1)^{k} \chi(\d_k^+X^H)
\end{eqnarray} 
Here we are employing  the pullback metric $F^\ast(g)$ in vicinity of $\d_1X \subset X$ to generate the volume form $d\mu_H$ on  $\d_1X^H$ and its normal curvature $K_H$. 
\end{cor}

\section{Probable  Refinements of the Equivariant Morse Formula}

One can refine the definition Birnside ring $A(G)$ in a number of ways. For example, following Dieck [D1, I.10.3], one can introduce the \emph{component category} $\Pi_0(G, X)$ of a $G$-space $X$ whose objects are $G$-maps $x: G/H \to X$. A morphism $\s$ from $x: G/H \to X$ to $y: G/K \to X$ is a $G$-map $\s: G/H \to G/K$ such that $x$ and $y\circ \s$ are $G$-homotopic. Denote by $Is\, \Pi_0(G, X)$ the set of isomorphism classes $[x]$ of objects $x: G/H \to X$. Let $\mathcal A^G(X) := \Z[Is\, \Pi_0(G, X)]$, where $\Z[S]$ stands for a free abelian group with the basis $S$. In fact, there is a bijection $$Is\, \Pi_0(G, X) \to \coprod_{(H) \in conj(G)}\; WH\setminus \pi_0(X^H).$$ This construction defines a covariant functor $\mathcal A^G(\sim)$ on the category of $G$-spaces $X$ with values in the category of abelian groups, a functor which is sensitive to the connected component structure of the sets $\{X^H\}$. We conjecture that all previous results  can be restated in terms of $\mathcal A^G(\sim)$ along the lines of [LR].
\bigskip

However, we would like to stress a different generalization of $A(G)$ and to speculate about the corresponding Morse Formulas. In this generalization one pays a close attention to the $H$-representations $\{\psi^H_\a\}_\a$ arising in the normal bundles  $\nu(X^H, X)$.\smallskip

Let $G$ be a compact Lie group. For each $(H) \in \Phi(G)$, fix a set $\Psi(H)$ of distinct isomorphism classes of representations $\{\psi^H_\a: H \to GL(V_\a) \}_\a$ with the property  $V_\a^H = \{0\}$. Moreover, for any  $\psi_\a \in \Psi(H)$ and $K \subset H$, the representation $\tilde{Res}_K(\psi_\a): K  \to GL(V_\a/ V_\a^K)$ is required to be in $\Psi(K)$. We call a collection $\mathcal F := \{\psi_\a\}_{(H) \in \Phi(G),\; \a \in  \Psi(H)}$ of such representations a \emph{normal} family. 

Let $X^H_{\psi_a}$ denote the set of points in $X^H$ with the normal representations isomorphic to $\psi_\a$.

\begin{defn} Let $\mathcal F = \mathcal F(G)$ be a normal family of representations. Two compact smooth $G$-manifolds $X$ and $Y$ are $\mathcal F$-equivalent if, for each $(H) \in \Phi(G)$, and $\psi_a \in \Psi(H) \subset \mathcal(G)$,
$$ \chi(X^H_{\psi_a}) = \chi(Y^H_{\psi_a}).$$
We denote by $A(G,\; \mathcal F)$ the group of such equivalence classes. 
\end{defn}
\begin{conj} All the equivariant Morse formulas above are valid in the refined Burnside group $A(G,\; \mathcal F)$.
\end{conj}

\end{document}